\documentclass{article}

\usepackage{amsmath,amsthm,amssymb}
\usepackage{array,multirow}
\usepackage{bbm, enumerate}
\usepackage[subnum]{cases}
\usepackage[shortlabels]{enumitem}
\usepackage[yyyymmdd,hhmmss]{datetime}
\usepackage{hyperref}
\usepackage{tikz}
\usepackage{mathtools} 
\usetikzlibrary{arrows,shapes}
\usetikzlibrary{positioning}

\allowdisplaybreaks

\numberwithin{equation}{section}

\let\originalleft\left
\let\originalright\right
\renewcommand{\left}{\mathopen{}\mathclose\bgroup\originalleft}
\renewcommand{\right}{\aftergroup\egroup\originalright}

\newcommand{\n}[1]{\multicolumn{1}{|c|}{#1}}

\newtheorem{theo}{Theorem}[section]
\newtheorem{definition}[theo]{Definition}
\newtheorem{conjecture}[theo]{Conjecture}
\newtheorem{lemma}[theo]{Lemma}
\newtheorem{theorem}[theo]{Theorem}
\newtheorem{corollary}[theo]{Corollary}
\newtheorem{example}[theo]{Example}
\newtheorem{note}[theo]{Note}
\newtheorem{proposition}[theo]{Proposition}

\newcommand{\nW}[2]{\tikz[remember picture,anchor=base, baseline,inner xsep=0pt]{\node[box](#2){\vphantom{$\int$}$\underline{#1}$};}}
\newcommand{\nS}[1]{\tikz[remember picture,anchor=base, baseline,inner xsep=0pt]{\node[box](#1){\vphantom{$\int$}\textvisiblespace};}}

\pgfmathsetmacro{\myinnersep}{2}
\tikzset{
box/.style={
    inner sep=0,%
    outer sep=0,%
    minimum height=5mm,%
    align=center}
}

\begin{document}
\title{Expected value of letters of permutations with a given number of $k$-cycles}
\author{Peter Kagey}

\maketitle

\tableofcontents


In this paper, we study permutations $\pi \in S_n$ with exactly $m$ transpositions.
In particular, we are interested in the expected value of $\pi(1)$ when such
permutations are chosen uniformly at random. When $n$ is even, this expected value
is approximated closely by $(n+1)/2$, with an error term that is related to the number isometries
of the $(n/2-m)$-dimensional hypercube that move every face.
Furthermore, when $k \mid n$, this construction generalizes to allow us to compute
the expected value of $\pi(1)$ for permutations
with exactly $m$ $k$-cycles. In this case, the expected value has an
error term which is related instead to the number derangements of the
generalized symmetric group $S(k,n/k-m)$.

When $k$ does not divide $n$, the expected value of $\pi(1)$ is precisely
$(n+1)/2$.
Indirectly, this suggests the existence of a reversible algorithm
to insert a letter into a permutation which preserves the number of $k$-cycles,
which we construct.

\section{Background}
\label{section:background}

In 2010, Mark Conger \cite{conger} proved that a permutation
with $k$ descents has an expected first letter of $\pi(1) = k + 1$,
independent of $n$.
This paper has the same premise, but with a different permutation statistic:
the number of $k$-cycles of a permutation.

This section, (Section \ref{section:background}) provides an overview of where
we're headed, and includes an critical example that will hopefully spark the
reader's curiosity and motivate the remainder of the paper.

Section \ref{section:recursiveStructure} establishes some recurrence relations for
the number of permutations in $S_n$ with a given number of $k$-cycles.
It also contains a theorem that gives an explicit way to compute the expected
value of the first letter based on these counts.

Section \ref{section:wreathProduct} describes an explicit correspondence
between $k$-cycles of permutations in $S_{kn}$ and fixed points of elements of the generalized
symmetric group $(\mathbb{Z}/k\mathbb{Z}) \wr S_n$. Using generating functions and results
from the previous section, this shows that the expected value of $\pi(1)$ of a
permutation with a given number of $k$-cycles is intimately connected to the
number of derangements of a generalized symmetric group.

While Section \ref{section:wreathProduct} emphasizes the case of $S_{kn}$,
Section \ref{section:bijection} looks at $S_N$ where $k \nmid N$. Here, the
expected value of $\pi(1)$ is simply $(N+1)/2$, which agrees with the expected
value of the first letter of a uniformly chosen $N$-letter permutation with
no additional restrictions. This fact together with the main
theorem from Section \ref{section:recursiveStructure} implies the existence of a
bijection $\varphi_k \colon S_{N-1} \times [N] \rightarrow S_N$ that preserves the
number of $k$-cycles whenever $k \nmid N$. Section \ref{section:bijection}
constructs such a bijection explicitly, and proves that it has the desired
properties.

\subsection{Motivating Examples}

In support of the first examples, we start by defining the first bit of notation.
\begin{definition}
  Let $C_k(n,m)$ denote the number of permutations $\pi \in S_n$ such that
  $\pi$ has exactly $m$ $k$-cycles.
\end{definition}

These theorems---and many of the following lemmas---were discovered by looking
at examples such as the following, written in both one-line and cycle notation:

\begin{example}
\label{ex:fourLettersTwoCycles}
  There are $C_2(4,0) = 15$ permutations in $S_4$ with no $2$-cycles:
  \begin{alignat*}{5}
    1234 &= (1)(2)(3)(4) \hspace{0.5 cm}
    && 2314 = (312)(4)\hspace{0.5 cm}
    && 3124 = (321)(4)\hspace{0.5 cm}
    && 4123 = (4321) \\
    1342 &= (1)(423) \hspace{0.5 cm}
    && 2341 = (4123) \hspace{0.5 cm}
    && 3142 = (4213) \hspace{0.5 cm}
    && 4132 = (421)(3)\\
    1423 &= (1)(432) \hspace{0.5 cm}
    && 2413 = (4312) \hspace{0.5 cm}
    && 3241 = (2)(413) \hspace{0.5 cm}
    && 4213 = (2)(431) \\
    & \
    && 2431 = (3)(412)  \hspace{0.5 cm}
    && 3421 = (4132) \hspace{0.5 cm}
    && 4312 = (4231)
  \end{alignat*}
  There are $C_2(4,1) = 6$ permutations in $S_4$ with exactly one $2$-cycle:
  \begin{alignat*}{3}
    1243 &= (1)(2)(43) \hspace{1cm}
    && 2134 = (21)(3)(4) \\
    1324 &= (1)(32)(4) \hspace{1cm}
    && 3214 = (2)(31)(4)\\
    1432 &= (1)(3)(42) \hspace{1cm}
    && 4231 = (2)(3)(41)
  \end{alignat*}
  And there are $C_2(4,2) = 3$ permutations in $S_4$ with exactly two $2$-cycles,
  \begin{equation*}
    2143 = (21)(43) \hspace{1cm}
    3412 = (31)(42) \hspace{1cm}
    4321 = (32)(41).
  \end{equation*}
  By averaging the first letter over these examples, we can
  compute that \begin{alignat*}{3}
    &\mathbb{E}[\pi(1)\, |\, \pi \in S_4 \text{ has no } 2 \text{-cycles}\,]
      &&= \frac{3(1) + 4(2 + 3 + 4)}{15}
      &&= \frac{13}{5},
    \\
    &\mathbb{E}[\pi(1)\, |\, \pi \in S_4 \text{ has exactly } 1\ 2 \text{-cycle}]
      &&= \frac{3(1) + (2 + 3 + 4)}{6}
      &&= 2, \text{ and}
    \\
    &\mathbb{E}[\pi(1)\, |\, \pi \in S_4 \text{ has exactly } 2\ 2 \text{-cycles}\,]
      &&= \frac{2 + 3 + 4}{3}
      &&= 3.
  \end{alignat*}
\end{example}
The table in Figure \ref{fig:twoCyclesTable} gives the expected value of $\pi(1)$ given that
$\pi \in S_n$ and has exactly $m$ $2$-cycles in its cycle decomposition.
Notice that when $i$ is odd, row $i$ has a constant value of $(i+1)/2$. Also
notice that the number in position $(i,j)$ has the same denominator as the
number in position $(i+2, j+1)$, and that these denominators increase with $n$.
The sequence of denominators begins \begin{equation}
  1, 5, 29, 233, 2329, 27949, \dots,
\end{equation}
which agrees with the type B derangement numbers,
sequence A000354 in the On-Line Encyclopedia of Integer Sequences (OEIS) \cite{oeis}.
In other words, the denominators in the table appear to be related to
the symmetries of the hypercube that move every facet.
\begin{figure}[!ht]
  \[
  \begin{array}{|ll|l|l|l|l|l|l|l|l|l|l|}
  \hline
  & & \multicolumn{7}{|c|}{m} \\ \cline{3-9}
  & & 0 & 1 & 2 & 3 & 4 & 5 & 6 \\ \hline
  \multirow{8}{*}{$n$}
  & \n{1}  & 1/1          &            &          &        &      &     &     \\ \cline{2-4}
  & \n{2}  & 1/1          & 2/1        &          &        &      &     &     \\ \cline{2-4}
  & \n{3}  & 2/1          & 2/1        &          &        &      &     &     \\ \cline{2-5}
  & \n{4}  & 13/5         & 2/1        & 3/1      &        &      &     &     \\ \cline{2-5}
  & \n{5}  & 3/1          & 3/1        & 3/1      &        &      &     &     \\ \cline{2-6}
  & \n{6}  & 101/29       & 18/5       & 3/1      & 4/1    &      &     &     \\ \cline{2-6}
  & \n{7}  & 4/1          & 4/1        & 4/1      & 4/1    &      &     &     \\ \cline{2-7}
  & \n{8}  & 1049/233     & 130/29     & 23/5     & 4/1    & 5/1  &     &     \\ \cline{2-7}
  & \n{9}  & 5/1          & 5/1        & 5/1      & 5/1    & 5/1  &     &     \\ \cline{2-8}
  & \n{10} & 12809/2329   & 1282/233   & 159/29   & 28/5   & 5/1  & 6/1 &     \\ \cline{2-8}
  & \n{11} & 6/1          & 6/1        & 6/1      & 6/1    & 6/1  & 6/1 &     \\ \cline{2-9}
  & \n{12} & 181669/27949 & 15138/2329 & 1515/233 & 188/29 & 33/5 & 6/1 & 7/1 \\ \cline{2-9}
  & \n{13} & 7/1          & 7/1        & 7/1      & 7/1    & 7/1  & 7/1 & 7/1 \\ \hline
  \end{array}
  \]
  \caption{
    A table of the expected value of the first letter of $\pi \in S_n$ with
    exactly $m$ $2$-cycles,
    ${\mathbb{E}[\pi(1)\, |\, \pi \in S_n \text{ has exactly } m\ 2 \text{-cycles}\,]}$.
  }
  \label{fig:twoCyclesTable}
\end{figure}
\section{Structure of permutations with \texorpdfstring{$m$}{m} \texorpdfstring{$k$}{k}-cycles}
\label{section:recursiveStructure}
This section is about connecting the number of permutations with a given number of
$k$-cycles to the expected value of the first letter. Saying this, it is
appropriate to start with a 1944 theorem of Goncharov that, by the principle of
inclusion/exclusion, gives an explicit formula that counts the number of such
permutations.
\subsection{Counting permutations based on cycles}
\begin{theorem}[\cite{goncharov}, \cite{arratia}]
  \label{GoncharovTheorem}
  The number of permutations in $S_n$ with exactly $m$ $k$-cycles is given by
  the following sum, via the principle inclusion/exclusion:
  \begin{equation}
    C_k(n,m)
    = \frac{n!}{m!k^m}\sum_{i=0}^{\lfloor n/k \rfloor - m} \frac{(-1)^i}{i!\,k^i}.
  \end{equation}
  \label{eq:Goncharov}
\end{theorem}

\begin{corollary}
  For $k \nmid n$, there are exactly $n$ times as many permutations in $S_n$
  with exactly $m$ $k$-cycles than there are in $S_{n-1}$.
  When $k \mid n$, there is an explicit formula for the difference.
  \label{cor:bijectionDifference}
  \begin{numcases}{C_k(n, m) - nC_k(n-1,m) = }
    0 & $k \nmid n$ \label{eq:bijectionDifference0}
    \\
    \displaystyle\frac{n!(-1)^{\frac nk - m}}{(n/k)! \, k^{\frac nk}}\binom{n/k}{m} & $k \mid n$
    \label{eq:bijectionDifference1}
  \end{numcases}
  \label{eq:bijectionDifference}
\end{corollary}
\begin{proof}
  When $k \nmid n$,
  $\displaystyle \left\lfloor \frac{n}{k}\right\rfloor = \left\lfloor \frac{n-1}{k}\right\rfloor$,
  so the bounds on the sums are identical and the result follows directly \begin{align}
    \frac{n!}{m!k^m}\sum_{i=0}^{\lfloor n/k \rfloor - m} \frac{(-1)^i}{i!\,k^i}
    - n\left(\frac{(n-1)!}{m!k^m}\sum_{i=0}^{\lfloor (n-1)/k \rfloor - m} \frac{(-1)^i}{i!\,k^i}\right)
    = 0.
  \end{align}
  Otherwise, when $k \mid n$,
  $\displaystyle \left\lfloor \frac{n-1}{k}\right\rfloor = \frac{n}{k} - 1$, so
  \begin{align}
    \nonumber
    &\frac{n!}{m!k^m}\sum_{i=0}^{n/k - m} \frac{(-1)^i}{i!\,k^i}
    - n\left(\frac{(n-1)!}{m!k^m}\sum_{i=0}^{n/k - 1 - m} \frac{(-1)^i}{i!\,k^i}\right) \\[10pt]
    \nonumber
    &\hspace{4cm}= \frac{n!}{m!k^m}\left(\frac{(-1)^{n/k - m}}{(n/k - m)!k^{n/k - m}}\right) \\[10pt]
    \nonumber
    &\hspace{4cm}= \frac{n!(-1)^{n/k - m}}{(n/k - m)!m!k^{n/k}} \\[10pt]
    &\hspace{4cm}= \frac{n!(-1)^{\frac nk - m}}{(n/k)! \, k^{n/k}}\binom{n/k}{m}.
  \end{align}
\end{proof}
See Section \ref{section:bijection} for a bijective proof of Equation
\ref{eq:bijectionDifference0}.

\subsection{Permutations by first letter}
In order to compute the expected value of the first letter of a permutation,
it is useful to be able to compute the number of permutations that have a given
number of $k$-cycles \textit{and} a given first letter.
\begin{definition}
  Let $C_k^{(a)}(n,m)$ be the number of permutations $\pi \in S_n$ that
  have exactly $m$ $k$-cycles and $\pi(1) = a$.
\end{definition}
The expected value of $\pi(1)$ with a given number of
$k$-cycles
is
  \begin{equation}
    \mathbb{E}[\pi(1)\, |\, \pi \in S_n \text{ has exactly } m\ k \text{-cycles}\,] =
    \frac{1}{C_k(n, m)}\sum_{a=1}^n a C_k^{(a)}(n, m).
    \label{eq:weightedAverage}
  \end{equation}
The following three lemmas compute $C_k^{(a)}(n,m)$ from $C_k(n,m)$.
\begin{proposition}
  \label{cyc1Recurrence}
  For all $k > 1$, the number of permutations in $S_n$ starting with $1$ and
  having $m$ $k$-cycles is equal to the number of permutations in $S_{n-1}$ with
  $m$ $k$-cycles: \begin{equation}
    C_k^{(1)}(n,m) = C_k(n-1, m).
  \end{equation}
\end{proposition}
\begin{proof}
  The straightforward bijection from $\{\pi \in S_n : \pi(1) = 1\}$ to $S_{n-1}$
  given by deleting $1$ and relabeling preserves the number of $k$-cycles for
  $k > 1$.
\end{proof}
\begin{proposition}
  \label{allSame}
  For all $a, b \geq 2$, the number of permutations having $k$-cycles
  and starting with $a$ are the same as the number of those starting with $b$:
  \begin{equation}
    C_k^{(2)}(n,m) = \cdots = C_k^{(a)}(n,m) = \cdots = C_k^{(b)}(n,m) = \cdots = C_k^{(n)}(n,m).
  \end{equation}
\end{proposition}
\begin{proof}
  Since the permutations under consideration do not fix $1$,
  conjugation by $(ab)$ is an isomorphism which takes all words starting
  with $a$ to words starting with $b$ without changing the cycle structure.
\end{proof}
\begin{lemma}
  \label{cycRecurrenceWithFixedBeginning}
  For all $2 \leq a \leq n$, \begin{align}
    C_k^{(a)}(n,m) = \frac{C_k(n, m) - C_k(n-1, m)}{n - 1}.
  \end{align}
\end{lemma}
\begin{proof}
  Since \begin{equation}
    C_k(n, m) = C_k^{(1)}(n, m) + C_k^{(2)}(n, m) + \dots + C_k^{(n)}(n, m),
  \end{equation} using Proposition \ref{allSame}, for the last $(n-1)$ terms,
  this can be rewritten as \begin{equation}
    C_k(n, m) = C_k^{(1)}(n, m) + (n-1)C_k^{(a)}(n, m).
  \end{equation}
  Solving for $C_k^{(a)}(n, m)$ and using the substitution from Proposition
  \ref{cyc1Recurrence} gives the desired result.
\end{proof}
Now, equipped with explicit formulas for $C_k^{(a)}(n,m)$ and $C_k(n,m)$, we
can compute the expected value of $\pi(1)$ for $\pi \in S_n$ with exactly $m$
$k$-cycles.
\subsection{Expected value of first letter}
\begin{theorem}
  \label{firstTheorem}
  For $k > 1$, the expected value of the first letter of a permutation
  $\pi \in S_n$ with $m$ $k$-cycles is given by \begin{align}
    &\mathbb E[\pi(1)\ \mid\ \pi \in S_n \text{ has exactly } m\ k\text{-cycles}\,]
    \nonumber
    \\
    & \hspace{4cm} = \frac n2\left(1 - \frac{C_k(n-1,m)}{C_k(n,m)}\right) + 1.
  \end{align}
\end{theorem}
\begin{proof}
  Using Proposition \ref{allSame}, we can consolidate all but the first term of
  the sum in Equation \ref{eq:weightedAverage} \begin{align}
    &\sum_{a = 1}^n aC_k^{(a)}(n, m) \\
    &\hspace{1cm}= C_k^{(1)}(n,m) + \sum_{a = 2}^n aC_k^{(n)}(n, m) \\
    &\hspace{1cm}= C_k^{(1)}(n,m) + \frac{(n-1)(n+2)}{2} C_k^{(n)}(n, m) \\
    &\hspace{1cm}=
    C_k(n-1,m) + \frac{(n-1)(n+2)}{2}\left(
      \frac{C_k(n, m) - C_k(n-1, m)}{n - 1}
    \right) \\
    &\hspace{1cm}= \left(\frac{n}{2} + 1\right) C_k(n,m) - \frac n2C_k(n-1,m).
  \end{align}
  Dividing by $C_k(n,m)$ yields the result.
\end{proof}
\begin{corollary}
  \label{cor:kNotDivideN}
  When $k \nmid n$, $C_k(n,m) = nC_k(n-1,m)$ by Equation \ref{eq:bijectionDifference0},
  so \begin{equation}
    \mathbb E[\pi(1)\ \mid\ \pi \in S_n \text{ has exactly } m\ k\text{-cycles}\,] = \frac{n}{2}\left(1 - \frac{1}{n}\right) + 1 = \frac{n+1}{2}.
  \end{equation}
\end{corollary}
Together with Theorem \ref{GoncharovTheorem}, this theorem and its corollary
provides our first formula for the expected value of $\pi(1)$ that
performs exponentially better than brute force.

\subsection{Identities for counting permutations with given cycle conditions}
Both in practical terms (if computing the expected value of $\pi(1)$ by hand or
optimizing an algorithm) and in a theoretical sense, the following recurrence is
simple and useful.

\begin{lemma}
  \label{cycleRecursion}
  For $n < mk$ or $m < 0$, $C_k(n, m) = 0$. Otherwise,
  for all $k, m \geq 1$ \begin{equation}
    mC_k(n, m) = (k-1)!\binom{n}{k}C_k(n-k, m-1).
  \end{equation}
\end{lemma}
While this can be proven directly by the algebraic manipulation of the
identity in Theorem \ref{GoncharovTheorem}, a bijective proof has been
included here because it is natural and may be of interest.
\begin{proof}
  Let \begin{equation}
  \mathcal C_k(n, m) = \{ \pi \in S_n\,\mid\,\pi \text{ has exactly } m\ k \text{-cycles}\}.
  \end{equation}
  Then consider the two sets, whose cardinalities match the left- and
  right-hand sides of the equation above:
  \begin{align}
    X^{L}_{n,m,k} &= \{ (\pi, c) \mid \pi \in \mathcal C_k(n, m), c \text{ a distinguished } k\text{-cycle of } \pi \}. \\
    X^{R}_{n,m,k} &= \{ (\sigma, d) \mid \pi \in \mathcal C_k(n-k, m-1), d \text{ an } n\text{-ary necklace of length } k\}.
  \end{align}
  The first set, $X^{L}_{n,m,k}$, is constructed by taking a permutation in
  $\mathcal C_k(n,m)$ and choosing one of its $m$ $k$-cycles to be distinguished, so
  $\#X^{L}_{n,m,k} = mC_k(n,m)$.

  In the second set, $X^{R}_{n,m,k}$, the two parts of the tuple are independent.
  There are $C_k(n-k, m-1)$ choices for the permutation $\sigma$ and $(k-1)!\binom{n}{k}$
  choices for the necklace $d$.
  Thus $\#X^{R}_{n,m,k} = (k-1)!\binom{n}{k}C_k(n-k, m-1)$.

  Now, consider the map $\varphi \colon X^{L}_{n,m,k} \rightarrow X^{R}_{n,m,k}$
  which removes the distinguished $k$-cycle and relabels the remaining $n - k$
  letters as $\{1, 2, \dots, n - k\}$, preserving the relative order:
  \begin{equation}
    (\pi_1\pi_2 \cdots \pi_\ell, \pi_i) \xmapsto{\varphi} (\pi'_1\pi'_2 \cdots \pi'_{i-1}\pi'_{i+1} \cdots \pi'_\ell, \pi_i)
  \end{equation} where $\pi'_i$ is $\pi_i$ after relabeling.

  By construction, $\sigma$ has one fewer $k$-cycle and $k$ fewer letters
  than $\pi$.

  The inverse map is similar. To recover $\pi$, increment the letters of $\sigma$ appropriately
  and add the necklace $d$ back in as the distinguished cycle.
  Thus $\varphi$ is a bijection and $\#X^{L}_{n,m,k} = \#X^{R}_{n,m,k}$.
\end{proof}
\begin{example}
  Suppose $\pi = (423)\mathbf{(61)}(75)$ in cycle notation with $(61)$ distinguished.
  Then \begin{align}
    \varphi((423)(61)(75), (61)) = ((312)(54), (61))
  \end{align} under the bijection $\varphi$, described in the proof of
  Lemma \ref{cycleRecursion}.
\end{example}
The recurrence in Lemma \ref{cycleRecursion} suggests that understanding
$C_k(n,m)$ is related to understanding $C_k(n-km, 0)$,
the permutations of $S_{n-km}$ with no $k$-cycles.
On the other hand,
Corollary \ref{cor:bijectionDifference} suggests that the case where $k \mid n$
has some of the most intricate structure.
We can, of course, combine these two observations and analyze the case of
$C_k(kn, 0)$, which has a particularly simple generating function, which will
show up again in a different guise.
\begin{lemma}
  For $k \geq 2$,
  \label{cycleEGF}
  \begin{equation}
    \sum_{n=0}^\infty \frac{C_k(kn, 0)k^n}{(kn)!}x^n
    = \frac{\exp(-x)}{1-kx}.
  \end{equation}
\end{lemma}
\begin{proof}
  By substitution of $C_k(kn, 0)$ via the identity in Theorem \ref{GoncharovTheorem},
  \begin{align}
    \sum_{n=0}^\infty \frac{C_k(kn, 0)k^n}{(kn)!}x^n
    &= \sum_{n=0}^\infty \sum_{i=0}^n \frac{(-1)^i}{k^i i!}k^nx^n \\
    &= \sum_{n=0}^\infty \sum_{i=0}^n \frac{(-x)^i}{i!}(kx)^{n-i} \\
    &= \left(\sum_{n=0}^\infty \frac{(-x)^n}{n!}\right) \left(\sum_{n=0}^\infty (kx)^n\right) \\
    &= \frac{\exp(-x)}{1-kx}.
  \end{align}
\end{proof}

This section allowed for the practical computation of the expected value of
$\pi(1)$ with a given number of $k$-cycles, but leaves the observation about
Figure \ref{fig:twoCyclesTable} unexplained. The following section will
explain the connection between the expected values of $\pi(1)$ and the
facet-derangements of the hypercube.

\section{Connection with the generalized symmetric group}
\label{section:wreathProduct}
This section explains the connection between the expected value of $\pi(1)$
given that $\pi$ has exactly $m$ $2$-cycles and the facet-derangements of the
hypercube, by telling the more general story of derangements of the generalized
symmetric group. Thus it is appropriate to start this section by defining
both the generalized symmetric group and its derangements.

\subsection{Derangements of the generalized symmetric group}
\begin{definition}
  The \textbf{generalized symmetric group} $S(k,n)$ is the wreath product
  $(\mathbb{Z}/k\mathbb{Z}) \wr S_n$, which in turn is a semidirect product
  $(\mathbb{Z}/k\mathbb{Z})^n \rtimes S_n$.
\end{definition}


A natural way of thinking about the symmetric group $S_n$ is by considering
how the elements act on length-$n$ sequences by permuting the indices.
Informally, we can think about the generalized symmetric group $S(k,n)$ in an
essentially similar way: each element consists of an ordered pair in
$(\mathbb{Z}/k\mathbb{Z})^n \rtimes S_n$, where $(\mathbb{Z}/k\mathbb{Z})^n$ gives information about
what to add componentwise, and $S_n$ gives information about how to rearrange
afterward.

\begin{example}
  Consider the generalized permutation \[
    (\underbrace{(1,3,0)}_{\in (\mathbb{Z}/4\mathbb{Z})^3}, \underbrace{(23)}_{\in S_3}) \in S(4,3).
  \]
  It acts on the sequence $(0,1,1) \in (\mathbb{Z}/2\mathbb{Z})^3$ first by adding
  element-wise, and then permuting: \begin{align}
    \underbrace{((1,3,0),(23))}_{\in S(k,n)} \cdot (0,1,1)
    = \underbrace{(23)}_{\in S_3} \cdot (1+0,3+1,0+1)
    = (23) \cdot (1,0,1)
    = (1,1,0).
  \end{align}
\end{example}

When $k = 1$, the sequence $(\mathbb{Z}/1\mathbb{Z})^n$ is trivially the zero
sequence, so $S(1,n) \cong S_n$.
When $k = 2$, $S(2,n)$ is the hyperoctahedral group that we brushed up against
in Figure \ref{fig:twoCyclesTable}:
the group of symmetries of the $n$-dimensional hypercube.
When $k \geq 3$, $S(k,n)$ does not have such an immediate geometric interpretation,
but it is precisely the right analog for the expected value of $\pi(1)$ when
$\pi$ has a given number of $k$-cycles.

\begin{definition}
  A \textbf{derangement} or \textbf{fixed-point-free element}
  of the generalized symmetric group is an element
  $((x_1,\dots,x_n),\pi) \in S(k,n)$ such that for all $i$,
  either $\pi(i) \neq i$ or $x_i \neq 0$.
\end{definition}

That is, when a derangement acts on a sequence in the manner described above,
it changes the position or the value of every term in the sequence.
When $k = 1$ and $S(1,n) \cong S_n$, this recovers the usual sense of a
derangement in $S_n$: a permutation with no fixed points.
In terms of the hyperoctahedral group, $S(2,n)$, a derangement is a symmetry of
the $n$-cube that moves each $(n-1)$-dimensional face.

\begin{example}
  The element $((1,3,0), (23)) \in S(4,3)$ is a derangement because it
  increments the first term and swaps the second and third terms---thus
  changing the position or value for each term.
\end{example}

The number of derangements of the generalized symmetric group can be described
by an explicit sum via the principle of inclusion/exclusion, and it has a
particularly elegant exponential generating function.

\begin{theorem}[\cite{assaf}] 
  \label{derangementEGF}
  For $k > 1$, the number of derangements of the generalized symmetric group $S(k,n)$ is
  \begin{equation}
    D(k,n) = k^n n!\sum_{i=0}^n \frac{(-1)^i}{k^i i!}.
  \end{equation} which has exponential generating function
  \begin{equation}
    \sum_{n=0}^\infty \frac{D(k,n)}{n!}x^n = \frac{\exp(-x)}{1 - kx}.
  \end{equation}
\end{theorem}
Notice that this agrees identically with the generating function in
Lemma \ref{cycleEGF}, which is our first hint in explaining the connection
between $k$-cycles in permutations and fixed points in elements of the
generalized symmetric group.
\subsection{Permutation cycles and derangements}
\begin{lemma}
  \label{CIdentity}
  For $k \geq 1$, the number of permutations with $kn + km$ letters and
  $m$ $k$-cycles is
  \begin{equation}
    C_k(k(n + m), m) = \binom{kn+km}{kn}C_k(kn,0)\frac{(km)!}{k^mm!}.
  \end{equation}
\end{lemma}
\begin{proof}[Algebraic proof]
  This will proceed by induction on $m$. The base case is clear when $m=0$,
  so suppose that the lemma is true up to $m-1$, that is \begin{align}
    C_k(k(n + m - 1), m - 1)
    &=
    \frac{(km-k)!}{k^{m-1}(m-1)!}\binom{kn + km - k}{kn}C_k(kn, 0). \\
    &= \frac{(kn + km - k)!}{k^{m-1}(m-1)!(kn)!}C_k(kn, 0).
    \label{eq:inductionStep1}
  \end{align}
  Rearranging Lemma \ref{cycleRecursion}, \begin{align}
    C_k(k(n + m), m)
    &= \frac{(k-1)!}{m}\binom{k(n + m)}{k}C_k(k(n+m-1), m-1) \\
    &= \frac{(kn+km)!}{km(kn+km-k)!}C_k(k(n+m-1), m-1).
    \label{eq:inductionStep2}
  \end{align}
  Now, notice there is a $(kn + km - k)!$ term in the numerator of Equation \ref{eq:inductionStep1} and
  the denominator of Equation \ref{eq:inductionStep2}, so substituting and simplifying yields \begin{equation}
    C_k(k(n + m), m)
    = \frac{(kn+km)!}{k^m m!(kn)!}C_k(kn, 0),
  \end{equation}
  as desired.
\end{proof}
\begin{proof}[Combinatorial proof]
  This lemma lends itself to a combinatorial proof.
  The left hand side of the equation counts the number of permutations in
  $S_{kn+km}$ with exactly $m$ $k$-cycles.
  The right hand side of the equation says that this is the number of ways to
  choose $kn$ letters in the permutation that will not be in $k$-cycles,
  and for each of these, there are $C_k(kn,0)$ ways to arrange these such that
  they have no $k$-cycles.
  This leaves over $km$ letters, of which there are
  $(km)!/(k^mm!)$ ways to write them as products of $m$
  disjoint $k$-cycles.
\end{proof}
The following lemma uses the above identities to establish that
the proportion of permutations in the symmetric group $S_{kn}$
with exactly $m$ $k$-cycles is equal to
the proportion of elements in the generalized symmetric group $S(k,n)$
with exactly $m$ fixed points.
\begin{lemma}
  \label{lem:fixedPointsAndCycles}
  For $k \geq 2$,
 \begin{equation}
    \frac{C_k(kn, m)}{(kn)!} = \binom nm\frac{D(k, n - m)}{k^nn!}.
    \label{eq:fixedPointsAndCycles}
  \end{equation}
\end{lemma}
\begin{proof}
  By solving for $D(k,n-m)$ on the right hand side and substituting $n + m$ for $n$, it is
  enough to show that the exponential generating function for $D(k,n)$
  (as shown in Theorem \ref{derangementEGF}) is also
  the exponential generating function for \begin{equation}
    C_k(kn+km, m) \frac{m!n!k^{n+m}}{(kn + km)!}.
  \end{equation}
  By the identity in Lemma \ref{CIdentity},
  \begin{align}
    &\sum_{n=0}^\infty C_k(kn+km, m) \frac{m!n!k^{n+m}}{(kn + km)!}\frac{x^n}{n!} \\
    &\hspace{2cm} = \sum_{n=0}^\infty \frac{(km)!}{m!k^m}\binom{kn + km}{kn}C_k(kn, 0) \frac{m!n!k^{n+m}}{(kn + km)!}\frac{x^n}{n!} \\
    &\hspace{2cm} = \sum_{n=0}^\infty C_k(kn, 0)\frac{k^nx^n}{(kn)!} \\
    &\hspace{2cm} = \frac{\exp(-x)}{1-kx},
  \end{align}
  with the final equality being the identity in Lemma \ref{cycleEGF}.
\end{proof}

\subsection{Expected value of letters of permutations}
We now have the ingredients we need to prove the pattern that we observed in
Figure \ref{fig:twoCyclesTable} that purported to show a relationship between
permutations given number of $2$-cycles and derangements of the hyperoctahedral
group. These ingredients come together in the following theorem, which
establishes the more general relationship between permutations with a given
number of $k$-cycles and derangements of the generalized symmetric
group, $S(k,n)$.

\begin{theorem}
  \label{secondTheorem}
  The expected value of the first letter of a permutation
  $\pi \in S_{kn}$ with exactly $m$ $k$-cycles, where $k > 1$ and $0 \leq m \leq n$, is
  \begin{equation}
    \label{eq:main}
    \mathbb{E}[\pi(1)\, |\, \pi \in S_{kn} \text{ has exactly } m\ k \text{-cycles}\,] =
    \frac{kn + 1}{2} + \frac{(-1)^{n-m}}{2 D(k, n - m)}
  \end{equation}
  where $D(k, n)$ is the number of derangements of the generalized symmetric
  group $S(k, n) = (\mathbb{Z}/m\mathbb{Z}) \wr S_n$.
\end{theorem}
\begin{proof}
  Inverting the identity in Lemma \ref{lem:fixedPointsAndCycles}, yields \begin{equation}
    \frac{\frac{(kn)!}{n!k^n}\binom n m}{C_k(kn, m)}
    = \frac{1}{D(k, n - m)}.
  \end{equation}
  Multiplying through by $(-1)^{n - m}$ to match the right hand side of Equation
  \ref{eq:main}, together with some small manipulations yields \begin{equation}
    1 - \frac{C_k(kn, m) - (-1)^{n - m}\frac{(kn)!}{n!k^n}\binom n m}{C_k(kn, m)}
    = \frac{(-1)^{n-m}}{D(k, n - m)}.
  \end{equation}
  Now adding $kn + 1$ and dividing by $2$ yields \begin{align}
    &\frac{kn}{2}\left(1 - \frac{C_k(kn, m) - (-1)^{n - m}\frac{(kn)!}{n!k^n}\binom n m}{knC_k(kn, m)}\right) + 1
    \nonumber
    \\
    & \hspace{2cm} = \frac{kn+1}{2} + \frac{(-1)^{n-m}}{2 D(k, n - m)},
  \end{align}
  which gives the right hand side as desired.
  Since the numerator on the left hand
  side is equal to ${kn C_k(kn - 1, m)}$ by Equation \ref{eq:bijectionDifference},
  the proof then follows from by Theorem \ref{firstTheorem}.
\end{proof}

With the expected value of the first letter found,
we can generalize this one more step to find the expected value of the
$i$-th letter of these permutations.

\begin{corollary}
  The expected value of the $i$-th letter of a permutation in $S_{kn}$ with
  exactly $m$ $k$-cycles, where
  $n \in \mathbb N_{>0}$,
  $k > 1$,
  $1 \leq i \leq kn$, and
  $0 \leq m \leq n$, is
  \[\mathbb{E}[\pi(i)\, |\, \pi \in S_{kn} \text{ has exactly } m\ k \text{-cycles}\,]
  = \frac{kn+1}{2}+\frac{(-1)^{n-m}}{2D(k,n-m)}\frac{kn+1-2i}{kn-1}.\]
\end{corollary}
\begin{proof}
  Denote by $N$ the number of permutations in $S_{kn}$
  with $m$ $k$-cycles where $1$ is a fixed point;
  denote by $M$ the number of permutations in $S_{kn}$
  with $m$ $k$-cycles where $\pi(1) = a \neq 1$.
  Note that while $N$ and $M$ implicitly depend on $m$, $n$, and $k$,
  $M$ does not depend on $a$ by Proposition \ref{allSame}.

  Thus \begin{align}
    \notag
    &\mathbb{E}[\pi(1)\, |\, \pi \in S_{kn} \text{ has exactly } m\ k \text{-cycles}\,] \\
    \notag
    &\hspace{1.5cm}
    = \frac{1}{N + (kn-1)M}\left(N + \sum_{a = 2}^{kn} aM \right) \\
    &\hspace{1.5cm}
    = \frac{1}{N + (kn-1)M}\left( N + \left(\frac{kn(kn+1)}{2} - 1\right)M \right).
  \end{align}
  More generally, if we conjugate with $(1i)$ then
  \begin{align}
    \notag
    &\mathbb{E}[\pi(i)\, |\, \pi \in S_{kn} \text{ has exactly } m\ k \text{-cycles}\,] \\
    \notag
    &\hspace{1.5cm} = \frac{1}{N + (kn-1)M}\left( N + \sum_{a \neq i} aM \right) \\
    \label{eq:affineFunction}
    &\hspace{1.5cm} = \frac{1}{N + (kn-1)M}\left(iN + \left(\frac{kn(kn+1)}{2} - i\right)M \right).
  \end{align}

  We can extend the function
  $\mathbb{E}[\pi(i)\, |\, \pi \in S_{kn} \text{ has exactly } m\ k \text{-cycles}\,]$
  to a function $f(n,k,m,i)$ where $i \in \mathbb Q$ is not necessarily an integer.
  As can be seen in Equation \ref{eq:affineFunction}, $f$ is affine function in $i$.
  By Theorem \ref{secondTheorem}, when $i = 1$, \[
    f(n,k,m,1) =
    \frac{kn + 1}{2} + \frac{(-1)^{n-m}}{2 D(k, n - m)}.
  \]
  When $i = (kn+1)/2$ yields \[
    f(n,k,m,(kn+1)/2) = \frac{kn+1}{2}.
  \]
  Because $f(n,k,m,i)$ is affine in $i$, it is enough to use linear
  interpolation and extrapolation to compute $f$ for arbitrary $i$.
  This can be done by scaling the
  $\displaystyle \frac{(-1)^{n-m}}{2 D(k, n - m)}$ term
  by an affine function of $i$ which is $1$ when $i=1$ and which vanishes
  when $i = (kn+1)/2$, namely $\displaystyle \frac{kn + 1 - 2i}{kn - 1}$, as
  desired.
\end{proof}

\begin{example}
  For $n = 2$, $k = 2$, and $m = 0$ the expected value of the first letter in a
  permutation in $S_{nk} = S_4$ with no $k=2$-cycles is $\displaystyle\frac{13}{5}$, as shown
  in Example \ref{ex:fourLettersTwoCycles}. This agrees with
  Theorem \ref{secondTheorem}:
  \begin{equation}
    \frac{kn + 1}{2} + \frac{(-1)^{n-m}}{2 D(k, n - m)} = \frac{4 + 1}{2} + \frac{(-1)^{2-0}}{2 D(2, 2 - 0)} = \frac{5}{2} + \frac{1}{10} = \frac{13}{5},
  \end{equation} since $D(2,2) = 5$ as illustrated in Figure \ref{fig:squareDerangements}.
  \begin{figure}[!ht]
    \centering
    \begin{tikzpicture}
      \draw (0,0) rectangle (2,2);
      \node[draw, fill=white, circle] at (1,2) {2}; 
      \node[draw, fill=white, circle] at (1,0) {4}; 
      \node[draw, fill=white, circle] at (0,1) {3}; 
      \node[draw, fill=white, circle] at (2,1) {1}; 
    \end{tikzpicture}
    ~
    \begin{tikzpicture}
      \draw (0,0) rectangle (2,2);
      \node[fill=white, circle] at (1,2) {3}; 
      \node[fill=white, circle] at (1,0) {1}; 
      \node[fill=white, circle] at (0,1) {4}; 
      \node[fill=white, circle] at (2,1) {2}; 
    \end{tikzpicture}
    ~
    \begin{tikzpicture}
      \draw (0,0) rectangle (2,2);
      \node[fill=white, circle] at (1,2) {4}; 
      \node[fill=white, circle] at (1,0) {2}; 
      \node[fill=white, circle] at (0,1) {1}; 
      \node[fill=white, circle] at (2,1) {3}; 
    \end{tikzpicture}
    ~
    \begin{tikzpicture}
      \draw (0,0) rectangle (2,2);
      \node[fill=white, circle] at (1,2) {1}; 
      \node[fill=white, circle] at (1,0) {3}; 
      \node[fill=white, circle] at (0,1) {2}; 
      \node[fill=white, circle] at (2,1) {4}; 
    \end{tikzpicture}

    \begin{tikzpicture}
      \draw (0,0) rectangle (2,2);
      \node[fill=white, circle] at (1,0) {2}; 
      \node[fill=white, circle] at (1,2) {4}; 
      \node[draw, fill=white, circle] at (0,1) {3}; 
      \node[draw, fill=white, circle] at (2,1) {1}; 

    \end{tikzpicture}
    ~
    \begin{tikzpicture}
      \draw (0,0) rectangle (2,2);
      \node[fill=white, circle] at (1,0) {3}; 
      \node[fill=white, circle] at (1,2) {1}; 
      \node[fill=white, circle] at (0,1) {4}; 
      \node[fill=white, circle] at (2,1) {2}; 
    \end{tikzpicture}
    ~
    \begin{tikzpicture}
      \draw (0,0) rectangle (2,2);
      \node[draw, fill=white, circle] at (1,0) {4}; 
      \node[draw, fill=white, circle] at (1,2) {2}; 
      \node[fill=white, circle] at (0,1) {1}; 
      \node[fill=white, circle] at (2,1) {3}; 
    \end{tikzpicture}
    ~
    \begin{tikzpicture}
      \draw (0,0) rectangle (2,2);
      \node[fill=white, circle] at (1,0) {1}; 
      \node[fill=white, circle] at (1,2) {3}; 
      \node[fill=white, circle] at (0,1) {2}; 
      \node[fill=white, circle] at (2,1) {4}; 
    \end{tikzpicture}
    \caption{The $2^{2}2! = 8$ symmetries of a square with fixed sides circled.
      The square ($2$-dimensional hypercube) has symmetry group
      $S(2,2) = (\mathbb{Z}/2\mathbb{Z}) \wr \mathbb S_2$ and $D(2,2) = 5$ of these
      symmetries are derangements, meaning that they do not fix any sides.}
    \label{fig:squareDerangements}
  \end{figure}
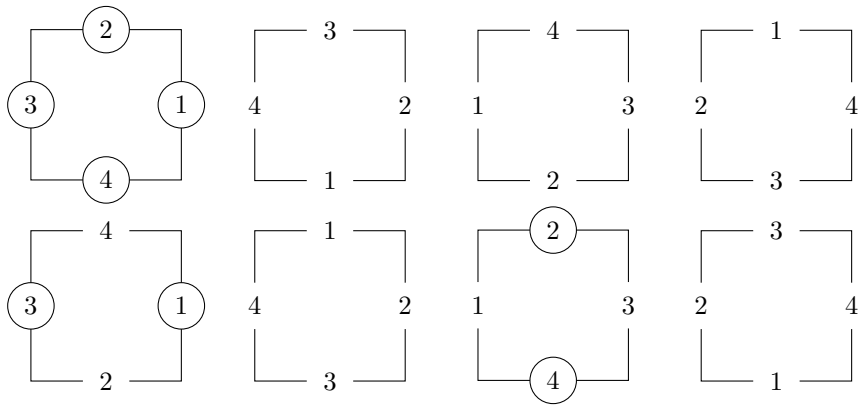
\end{example}

While Theorem \ref{firstTheorem} gave us our first way to efficiently compute
the expected value of the first letter of a permutation on $kn$ letters with
a given number of $k$-cycles, we can also compute this efficiently with
Theorem \ref{secondTheorem} by using the formulas for $D(k,n)$ in Theorem \ref{derangementEGF}.
But this is not the only reason that Theorem \ref{secondTheorem} is of interest;
because of the structure of the formula it provides, this theorem suggests other
quantitative and qualitative insights.

Recall that when there are no restrictions on a permutation $\pi \in S_{kn}$, the
first letter is equally likely to take on any value,
so $\mathbb{E}[\pi(1) \mid \pi \in S_{kn}] = (kn+1)/2$.
The first insight given by Theorem \ref{secondTheorem}
is that the expected value of $\pi(1)$ given some
number of $k$ cycles differs from $(kn+1)/2$ by at most $1/2$,
because $D(k,N) \geq 1$ for $k \geq 2$.
Secondly, since $D(k, N)$ increases as a function of $N$, the expected value
gets closer to $(kn+1)/2$ as the number of $k$-cycles decreases.
Lastly, the numerator of $(-1)^{n-m}$ in the second summand of
Equation \ref{eq:main} shows that the expected value of the first letter is
larger than $(kn+1)/2$ if and only if $n$ and $m$ have the same parity.

\section{A \texorpdfstring{$k$}{k}-cycle preserving bijection}
\label{section:bijection}
Motivated by Equation \ref{eq:bijectionDifference0}, this section describes a
family of bijections,
  \[\phi_k \colon S_{n-1} \times [n] \rightarrow S_n,\]
each of which preserves the number of $k$-cycles when $k \nmid n$.
Of course, there is no map that preserves the number of $k$-cycles when $k \mid n$. For example, a permutation in
$S_n$ consisting entirely of $k$-cycles contains $n/k$ $k$-cycles, while a permutation
in $S_{n-1}$ can contain at most $n/k - 1$ $k$-cycles by the pigeonhole principle.


Informally, these maps are defined by writing down a permutation $\sigma \in S_{n-1}$ in
\emph{canonical cycle notation}, incrementing all letters in $\sigma$
that are greater than or equal to $x \in [n]$,
inserting $x$ into the rightmost cycle, and then recursively moving letters into or out of subsequent cycles,
whenever a $k$-cycle is turned into a $(k+1)$-cycle or a $(k-1)$-cycle is
turned into a $k$-cycle.

\subsection{Example of recursive structure}

The definition of the map can look complicated, so it's worthwhile to start with
an example to give some sense of the overarching idea.

\begin{example}
  This example illustrates how the map $\phi_3$ inserts $I$ into
  the permutation $(D76)(E)(F32)(G91C)(K54)(LJ8)(MB)(NAH)$ while preserving the
  number of $3$-cycles.
  The maps $\phi_k$ and $\psi_k$ are the result of moving letters according
  to the arrows and are applied from right-to-left.
  (This example uses the convention that $1 < 2 < \dots < 9 < A < B < \dots < N$.)

  ~

  \begin{align*}
    &\phi_3(
      (D         7          6)
      (E                                          \nS{eEnd})
      (F         \nW{3}{3}  2                     \nS{2End})
      (G         \nW{9}{9}  1          \nW{C}{C})
      (K \nS{k5} 5          \nW{4}{4})
      (L \nS{lj} J          \nW{8}{8})
      (M \nS{mb} B                                \nS{bEnd})
      (N         \nW{A}{A}  H                     \nS{hEnd})
      ,
      \nW{I}{I}
    ) \\
    &\hspace{1cm} = (D76)(E3)(F29)(G1)(KC5)(L4J)(M8BA)(NHI)
  \end{align*}
  \begin{tikzpicture}[remember picture,overlay]
    \draw[thick, ->]  (I) to[out=90, in=90, looseness=3.4] node[midway, above]{$\phi_3$} (hEnd);
    \draw[thick, ->]  (A) to[out=90, in=90, looseness=2] node[midway, above]{$\phi_3$} (bEnd);
    \draw[thick, ->]  (8) to[out=90, in=90, looseness=2] node[midway, above]{$\psi_3$} (mb);
    \draw[thick, ->]  (4) to[out=90, in=90, looseness=2] node[midway, above]{$\psi_3$} (lj);
    \draw[thick, ->]  (C) to[out=90, in=90, looseness=2] node[midway, above]{$\psi_3$} (k5);
    \draw[thick, ->]  (9) to[out=90, in=90, looseness=2] node[midway, above]{$\phi_3$} (2End);
    \draw[thick, ->]  (3) to[out=90, in=90, looseness=2] node[midway, above]{$\phi_3$} (eEnd);
  \end{tikzpicture}

  \begin{align*}
    &\psi_3(
      (D         7         6)
      (E                   \nW{3}{3})
      (F \nS{f2} 2         \nW{9}{9})
      (G \nS{g1} 1                             \nS{1End})
      (K         \nW{C}{C} 5                   \nS{5End})
      (L         \nW{4}{4} J                   \nS{jEnd})
      (M         \nW{8}{8} B         \nW{A}{A})
      (N \nS{nh} H         \nW{I}{I})
    )\nS{outside}\\
    &\hspace{1cm} = ((D76)(E)(F32)(G91C)(K54)(LJ8)(MB)(NAH),I)
  \end{align*}
  \begin{tikzpicture}[remember picture,overlay]
    \draw[thick, ->]  (3) to[out=90, in=90, looseness=2] node[midway, above]{$\psi_3$} (f2);
    \draw[thick, ->]  (9) to[out=90, in=90, looseness=2] node[midway, above]{$\psi_3$} (g1);
    \draw[thick, ->]  (C) to[out=90, in=90, looseness=2] node[midway, above]{$\phi_3$} (1End);
    \draw[thick, ->]  (4) to[out=90, in=90, looseness=2] node[midway, above]{$\phi_3$} (5End);
    \draw[thick, ->]  (8) to[out=90, in=90, looseness=2] node[midway, above]{$\phi_3$} (jEnd);
    \draw[thick, ->]  (A) to[out=90, in=90, looseness=2] node[midway, above]{$\psi_3$} (nh);
    \draw[thick, ->]  (I) to[out=90, in=90, looseness=3.4] node[midway, above]{$\psi_3$} (outside);
  \end{tikzpicture}
\end{example}
Again, it is worth reemphasizing that the following definitions will follow the
convention that permutations are written in canonical cycle notation, \[
  \pi =
    \underbrace{(c^{(t)}_1\cdots c^{(t)}_{\ell_t})}_{c^{(t)}}
    \cdots
    \underbrace{(c^{(1)}_1\cdots c^{(1)}_{\ell_1})}_{c^{(1)}},
\] where cycle $c^{(i)} = (c^{(i)}_1 \cdots c^{(i)}_{\ell_i})$ has $\ell_i$ letters.
This means that the first letter in each cycle, $c^{(i)}_1$, is the largest letter in that cycle,
and that the cycles are ordered in increasing order by first letter when read from right-to-left:
$c^{(i+1)}_1 < c^{(i)}_1$ for all $i$.
\subsection{Formal definition and properties}
\begin{definition}
  Define
  $\phi_k \colon S_{n-1} \times [n] \mapsto S_{n}$ recursively as follows:
  \begin{equation}
    \phi_k(\emptyset, 1) = (1),
  \end{equation}
  and for $n>1$, $\pi \in S_{n-1}$, and $x \in [n]$,
  \begin{numcases}{\phi_k(\pi, x) =}
    \label{eq:phi_a}
    c^{(t)} \cdots c^{(1)}(x)                                                           & $x > c^{(1)}_1$            \\[10pt]
    \label{eq:phi_b}
    \phi_k(c^{(t)} \cdots c^{(2)}, c^{(1)}_2)(c^{(1)}_1 c^{(1)}_3 \cdots c^{(1)}_{k} x) & $\ell_1 = k$               \\[10pt]
    \label{eq:phi_c}
    \pi' (c^{(1)}_1 x' c^{(1)}_2\cdots c^{(1)}_{k-1} x)                                 & $\ell_1 = k-1, t > 1$ \\[10pt]
    \label{eq:phi_d}
    c^{(t)} \cdots c^{(2)} (c^{(1)}_1 \cdots c^{(1)}_{\ell_1} x)                        & otherwise.
  \end{numcases}
  Here, $\phi_k$ depends on the auxillary function
  $\psi_k \colon S_{n} \mapsto S_{n-1} \times [n]$,
  \begin{numcases}{\psi_k(\pi) =}
    \label{eq:psi_a}
    \left(c^{(t)} \dots c^{(2)}, c^{(1)}_1\right)                                                              & $\ell_1 = 1$             \\[10pt]
    \label{eq:psi_b}
    \left(\phi_k(c^{(t)} \cdots c^{(2)}, a_2^{(1)})(c^{(1)}_1 c^{(1)}_3 \dots c^{(1)}_k), c^{(1)}_{k+1}\right) & $\ell_1 = k + 1$         \\[10pt]
    \label{eq:psi_c}
    \left(\pi'(c^{(1)}_1 x' c^{(1)}_2 \dots c^{(1)}_{k-1}), c^{(1)}_k\right)                                   & $\ell_1 = k, t > 1$ \\[10pt]
    \label{eq:psi_d}
    \left(c^{(t)} \cdots c^{(2)}(c^{(1)}_1\cdots c^{(1)}_{\ell_1-1}),c^{(1)}_{\ell_1}\right)                   & otherwise,
  \end{numcases}
  and in both functions, $(\pi', x') = \psi(c^{(t)}\dots c^{(2)})$.
\end{definition}

\begin{note}
  Strictly speaking, $\phi_k$ and $\psi_k$ have an additional implicit
  parameter $n$, which indicates the size of permutation that these functions
  act on. Since the construction of these functions do not depend on $n$, this
  is suppressed in the notation.
\end{note}
The following theorem motivates this map,
and together with Lemma \ref{lem:bijection}, it implies Equation \ref{eq:bijectionDifference0}.
\begin{theorem}
  If $k \nmid n$, the number of $k$-cycles of $\pi \in S_{n-1}$ is equal to the number of $k$-cycles in
  $\phi_k(\pi,x)$.
\end{theorem}
\begin{proof}
  By construction, the maps $\phi_k$ and $\psi_k$ change the rightmost cycle
  into a (different) $k$-cycle if it was previously a $k$-cycle, and they
  change non-$k$-cycles into non-$k$-cycles, except for the case where there is
  one cycle remaining with length $k-1$ (in the case of $\phi$) or length $k$
  (in the case of $\psi$). These cases can only be achieved when $k \mid n$, by
  the following lemma.
\end{proof}
\begin{lemma}
  The number of letters in $\pi$ in (recursive) applications of $\phi_k$ and
  $\psi_k$ are of congruent to $n - 1 \bmod k$ and $n \bmod k$, respectively.
  Therefore, the only time that the input to $\phi_k$ can be a single cycle of
  length $k-1$ or the input to $\psi_k$ can be a single cycle of length $k$ is
  when $n \equiv 0\ (\bmod\ k)$.
\end{lemma}
\begin{proof}
  The proof proceeds by induction on the number of recursive iterations of
  $\phi_k$ and $\psi_k$. The base case is clear: on the first application of a
  map is always
  $\phi_k \colon S_{n-1} \times [n] \rightarrow S_n$, and the input permutation
  has $n-1$ letters by definition.

  Now, either we're finished, or we recurse
  (Equations \ref{eq:phi_b}, \ref{eq:phi_c}, \ref{eq:psi_b}, or \ref{eq:psi_c}),
  which we look at case-by-case.
  \begin{enumerate}[leftmargin=*, label={\textbf{Case \arabic*.}}]
    \item In Equation \ref{eq:phi_b}, the map $\phi_k$ sets aside $k$ letters from the input, so the number of letters in the recursive input to $\phi_k$ is also congruent to $n - 1 \bmod k$.
    \item In Equation \ref{eq:phi_c}, the map $\phi_k$ sets aside $k - 1$ letters from the leftmost cycle of the input. Since the number of letters in the original permutation was congruent to $n-1 \bmod k$, the number of letters in the permutation being input to $\psi_k$ is congruent to $n   \bmod k$.
    \item In Equation \ref{eq:psi_b}, the map $\psi_k$ sets aside $k + 1$ letters from the leftmost cycle of the input. Since the number of letters in the original permutation was congruent to $n   \bmod k$, the number of letters in the permutation being input to $\phi_k$ is congruent to $n-1 \bmod k$.
    \item In Equation \ref{eq:psi_c}, the map $\psi_k$ sets aside $k$ letters from the input, so the number of letters in the recursive input to $\psi_k$ is also congruent to $n \bmod k$.
  \end{enumerate}
\end{proof}
The following lemma provides a certain ``niceness'' property of the map,
which allows us to analyze it. In particular, all recursive inputs in both
$\phi_k$ and $\psi_k$ are written in canonical cycle notation.
\begin{lemma}
  The output of $\phi_k$ is in canonical cycle notation.
\end{lemma}
\begin{proof}
  Canonical cycle notation is preserved by construction.
  In particular, $\phi_k$ moves the first letter in any cycle, and
  Equation \ref{eq:phi_a} guards against inserting a number into a cycle that
  is bigger than the largest number already in the cycle.
  Similarly, $\psi_k$ only moves the first letter in the case of Equation
  \ref{eq:psi_a}, but in this case, the cycle only has one letter, so this is
  equivalent to deleting the cycle.
\end{proof}
\subsection{Inverting the bijection}
\begin{lemma}
  \label{lem:bijection}
  The maps
  $\phi_k \colon S_{n-1} \times [n] \rightarrow S_n$ and
  $\psi_k \colon S_n \rightarrow S_{n-1} \times [n]$ are inverse to one another.
\end{lemma}
\begin{proof}
  To prove this lemma, it suffices to show that $\psi_k \circ \phi_k = \operatorname{id}$
  by induction on the number of cycles of $\pi$. This will simultaneously prove
  that $\phi_k \circ \psi_k = \operatorname{id}$, because $S_{n-1} \times [n]$
  and $S_n$, both having $n!$ elements, have the same cardinality.

  When $\pi$ has no cycles, the base case is clear:
  $\psi_k(\phi_k(\emptyset, x)) = \psi_k((x)) = (\emptyset, x)$.

  Now there are five remaining cases to check, corresponding to each of the
  cases in the definition of $\phi_k(\pi, x)$
   \begin{enumerate}[leftmargin=*, label={\textbf{Case \arabic*.}}]
    \item Assume $x > c_1^{(1)}$, so that $\phi_k(\pi,x)$ is evaluated via Equation \ref{eq:phi_a}: \begin{align}
      \psi_k(\phi_k(\pi, x))
      &= \psi_k(c^{(t)}\cdots c^{(1)}(x)) \\
      &= (c^{(t)}\cdots c^{(1)},x) \\
      &= (\pi, x).
    \end{align}
    \item Assume $\ell_1 = k$, so that $\phi_k(\pi,x)$ is evaluated via Equation \ref{eq:phi_b}:\begin{align}
      \psi_k(\phi_k(\pi, x))
      &= \psi_k(\phi_k(c^{(t)} \cdots c^{(2)}, c^{(1)}_2)
      \underbrace{
        (c^{(1)}_1 c^{(1)}_3 \cdots c^{(1)}_{k} x)
      }_{\text{length } k}) \\
      &= (\pi'(c^{(1)}_1x'c^{(1)}_3 \dots c^{(1)}_k), x)
    \end{align}
    \item Assume $\ell_1 = k - 1$ and $t > 1$, so that $\phi_k(\pi,x)$ is evaluated via Equation \ref{eq:phi_c}:  \begin{align}
      \psi_k(\phi_k(\pi, x))
      &= \psi_k(\pi'
        \underbrace{
          (c^{(1)}_1 x' c^{(1)}_2\cdots c^{(1)}_{k-1} x)
        }_{\text{length } k + 1}
      )
    \end{align} where $(\pi', x') = \psi_k(c^{(t)} \dots c^{(2)})$.
    Therefore, this simplifies by Equation \ref{eq:psi_c}: \begin{align}
      \psi_k(\phi_k(\pi, x)) &= \left(
        \phi_k(\pi', x')
        (c^{(1)}_1 \cdots c^{(1)}_{k-1}),
        x
      \right) \\
      &= \Big(
        \underbrace{
          \phi_k(\psi_k(c^{(t)} \dots c^{(2)}))
        }_{c^{(t)} \dots c^{(2)}}
        \underbrace{
          (c^{(1)}_1 \cdots c^{(1)}_{k-1})
        }_{c^{(1)}},
        x
      \Big) \\
      &= (\pi, x),
    \end{align} because $\phi_k(\psi_k(c^{(t)} \dots c^{(2)})) = c^{(t)} \dots c^{(2)}$
    by the induction hypothesis on $t-1$ letters.
    \item Assume that $x > c_1^{(1)}$ and $\ell_1 \not\in \{k-1,k\}$, so that $\phi_k(\pi,x)$ is evaluated via Equation \ref{eq:phi_d}: \begin{align}
      \psi_k(\phi_k(\pi, x))
      &= \psi_k(c^{(t)} \cdots c^{(2)} (c^{(1)}_1 \cdots c^{(1)}_{\ell_1} x)) \\
      &= (c^{(t)}\cdots c^{(1)},x) \\
      &= (\pi, x).
    \end{align}
    \item Assume that $\ell_1 = k-1$ and $t = 1$, so that $\phi_k(\pi,x)$ is evaluated via Equation \ref{eq:phi_d}: \begin{align}
      \psi_k(\phi_k(\pi, x))
      &= \psi_k((c^{(1)}_1 \cdots c^{(1)}_{k-1  } x)) \\
      &= (c^{(1)},x) \\
      &= (\pi, x).
    \end{align}
  \end{enumerate}
\end{proof}
In this section we constructed a recursively-defined map and its inverse to
give a bijective proof that $C_k(n,m) = nC_k(n-1,m)$ when $k \nmid n$. This
is a novel, reversible algorithm for inserting a letters into a permutation
that preserves the number of $k$-cycles whenever possible.
\section{Further directions}
\label{section:furtherDirections}
In the introduction, we mentioned Conger's paper which analyzed how the number
of descents of a permutation affects the expected value of the first letter
of the permutation.
And similarly in the following sections, we looked at how the number of $k$-cycles
affects the expected value of the first letter of the permutation.
This section will principally look at the obvious generalization: given some
permutation statistic $\operatorname{stat}\colon S_n \rightarrow \mathbb Z$,
does the map \begin{equation}
  f(n,m) = \mathbb E[\pi(i) \mid \pi \in S_n, \operatorname{stat}(\pi)=m]
\end{equation} have any interesting structure?

But notice that the first letter of a permutation is itself a statistic, so
we can play a more general game. Given pairs of statistics
$(\operatorname{stat}_1, \operatorname{stat}_2)$, does the map
\begin{equation}
  g(n,m) = \mathbb E[\operatorname{stat}_1(\pi) \mid \pi \in S_n, \operatorname{stat}_2(\pi)=m]
\end{equation} have any interesting structure?

\subsection{FindStat database}
The result by Conger gives the expected value of $\pi(1)$ given
$\operatorname{des}(\pi)$, and this paper gave the expected value of
$\pi(1)$ given the number of $k$-cycles of $\pi$. Of course, it would be
interesting to do analogous analysis with other permutations. In particular,
the FindStat permutation statistics database \cite{FindStat} contains over
370 different permutation statistics, and many of these appear to have some
structure with respect to the expected value of the first letter of a
permutation.

\subsection{Mahonian statistics}
In particular, the family of Mahonian statistics may be fruitful to investigate.
Below, we have given conjectures about two: the major index and the inversion number.
Mahonian statistics are maps
$\operatorname{mah} \colon S_n \rightarrow \mathbb{N}_{\geq0}$ that are
equidistributed with the inversion number.\cite{Foata} That is, \[
  \#\{w \in S_n : \operatorname{mah}(w) = k\} =
  \#\{w \in S_n : \operatorname{inv}(w) = k\}.
\]
Naturally, all Mahonian statistics share the same generating function: \[
  \sum_{\sigma \in S_n} x^{\operatorname{mah}(\sigma)}
    = [n]_q!
    = \prod_{i=0}^{n-1}\sum_{j=0}^i(q^j).
\]

Because the expected value of the first letter is given by the weighted sum of
the permutations with $\operatorname{mah}(w) = k$ divided by the number of such
permutations, $\mathbb{E}[\pi(1)\, |\, \pi \in S_n, \operatorname{mah}(\pi) = k]$
has a denominator that is (a factor of) $M(n,k)$, the number of permutations
of $w \in S_n$ such that $\operatorname{inv}(w) = k$. For fixed $k$, these
satisfy a degree $k$ polynomial for all $n > k$. Notably, in the cases of
the major index and the inversion number, the numerators appear to satisfy
degree $k$ and degree $k-1$ polynomials respectively.

\begin{conjecture}
  For fixed $k$ and $n > k$, the expected value of the first letter of a
  permutation with a given number of inversions satisfies a rational function
  in $n$ given by \[
    \mathbb{E}[\pi(1)\, |\, \pi \in S_n, \operatorname{inv}(\pi) = k]
    = \frac{M(n+1,k)}{M(n,k)},
  \] where $M(n,k)$, as above, is the number of permutations $w \in S_n$ such
  that $\operatorname{inv}(w) = k$.
\end{conjecture}

\begin{conjecture}
  For fixed $k > 0$ and $n \geq k$,
  $\mathbb{E}[\pi(1)\, |\, \pi \in S_n, \operatorname{maj}(\pi) = k]$
  satisfies a rational function in $n$ that is $1/(k+1)$ times the quotient of a monic
  degree-$(k+1)$ polynomial by a monic degree-$k$ polynomial. Specifically,

  \begin{align}
    \mathbb{E}[\pi(1)\, |\, \pi \in S_n, \operatorname{maj}(\pi) = 1] &= \frac{1}{2}\left(\frac{n^2 + n - 2}{n-1}\right),
    \\[2mm]
    \mathbb{E}[\pi(1)\, |\, \pi \in S_n, \operatorname{maj}(\pi) = 2] &= \frac{1}{3}\left(\frac{n^3 - n - 6}{n^2 - n - 2}\right),
    \\[2mm]
    \mathbb{E}[\pi(1)\, |\, \pi \in S_n, \operatorname{maj}(\pi) = 3] &= \frac{1}{4}\left(\frac{n^4 + 6 n^3 - 13 n^2 - 18 n}{n^3 - 7n}\right),
    \text{ and}
    \\[2mm]
    \mathbb{E}[\pi(1)\, |\, \pi \in S_n, \operatorname{maj}(\pi) = 4] &= \frac{1}{5}\left(\frac{n^5 + 20 n^4 - 45 n^3 - 80 n^2 - 16 n}{n^4 + 2 n^3 - 13 n^2 - 14 n}\right).
  \end{align}
  Note that the denominator is given by an integer multiple of $M(n,k)$,
  a degree $k$ polynomial.
\end{conjecture}

\subsection{An elusive bijection}

  Let $F_k(n, m)$ be the number of elements of the generalized symmetric group
  $S(k,n) = (\mathbb{Z}/k\mathbb{Z}) \wr S_n$ with $m$ fixed points,
  and recall that $C_k(n,m)$ is the number of elements of $S_{kn}$ with $m$ $k$-cycles.
  Then for each pair of nonnegative integers $(\alpha, \beta)$
  with $\alpha, \beta \leq n$, then as Lemma \ref{lem:fixedPointsAndCycles} suggests,
  there exists a bijection of sets \begin{equation}
    C_k(n, \alpha) \times F_k(n, \beta) \rightarrow C_k(n, \beta) \times F_k(n, \alpha).
  \end{equation}
  This bijection has proven to be elusive to construct outside of the special
  cases where $n=1$ or $k=1$.
  Note that, the map cannot be a group automorphism of $S_{kn} \times S(k,n)$,
  because the identity of this group is in $C_k(n,0) \times F_k(n,n)$, so it
  cannot be preserved under this map.

  It would be especially interesting if there's a way to use the embedding of
  $(\mathbb{Z}/k\mathbb{Z}) \wr S_n$ into $S_{kn}$ as the centralizer of an
  element that is the product of $n$ disjoint $k$ cycles.

\section{Acknowledgments}
A special thanks to my advisor, Sami Assaf, for sharing the spark that she
found for this questions in a remark by Jim Pitman, and for her patient
guidance. This paper benefitted from the feedback from my colleague, Sam Armon
and his generosity, kindness, and sharp eye.
It is unlikely that this paper would have been written if not for the On-Line
Encyclopedia of Integer Sequences, which gave a several crucial hints,
especially around the pattern in Figure \ref{fig:twoCyclesTable}.
\bibliography{permutation_statistics_arxiv}
\end{document}